# Participation Factor-Based Adaptive Model Reduction for Fast Power System Simulation

Mahsa Sajjadi, *Student Member*, *IEEE*, Kaiyang Huang, *Student Member*, *IEEE*, Kai Sun, *Senior Member*, *IEEE*

Department of Electrical and Computer Engineering, University of Tennessee, Knoxville, TN, USA

msajjad1@vols.utk.edu, khuang12@vols.utk.edu, kaisun@utk.edu

*Abstract*—This paper describes an adaptive method to reduce a nonlinear power system model for fast and accurate transient stability simulation. It presents an approach to analyze and rank participation factors of each system state variable into dominant system modes excited by a disturbance so as to determine which regions or generators can be reduced without impacting the accuracy of simulation for a study area. In this approach, the generator models located in an external area with large participation factors are nonlinearly reduced and the rest of the generators will be linearized. The simulation results confirm that the assessment of the level of interaction between generators and system modes by participation factors is effective in enhancing the accuracy and speed of power system models. The proposed method is applied to the Northeastern Power Coordinating Council region system with 48-machine, 140-bus power system model and the results are compared with two cases including fully linearized model reduction and model reduction using the rotor angle deviation criteria.

## I. Introduction

Power system planners and utility operators rely on dynamic simulations to assess the dynamic behavior of a power system subject to a contingency and maintain a reliable and secure operating condition. The growth of massive electrical networks has necessitated more dynamic model analysis to timely diagnose eventual instability issues and neutralize the impact of such potential system instabilities following each disturbance.

Applying model reduction to a complex power system is one technique to improve the speed of online simulation. The frequently used solution for power system model reduction is to partition the whole network into two sections including 1) study area defined as to be the main goal of dynamic simulation study, in which all the details and nonlinear models are retained and 2) external area, where models can be truncated such that they reflect the overall performance of the remaining elements of the power system grid. Several tie-lines connect each part of the study area to the external area. Each of these tie-lines acts as a simple fictitious generator that represents the voltage magnitude and voltage angle of boundary bus between these two areas [1].

Several methods have so far been proposed for power system model reduction, some of which are: coherency-based approaches [2, 3], linear model reduction approaches such as low-rank Choleski factor method [4], dominating pole method [5], Krylov subspace and balanced truncation methods [6-8]. Each of these methods has its own unique features and limitations. These methods sometimes work well especially in small size disturbances, but they often work poorly, partly because of large and undesirable error value when a large disturbance happens in a real power system grid. To address this issue, a nonlinear model reduction can be used when the disturbance is large [9-11]. However, in this case, the speed of simulation will be significantly decreased compared to simulations using a linearly reduced model, which is an unwanted side effect.

In paper [9] and its following works [10, 11], an adaptive nonlinear model reduction method is proposed that requires knowing a good threshold ahead of time for simulated state variables, e.g. rotor angle derivations, to switch between linear and nonlinear model reduction algorithms. However, such a threshold is often system-specific and needs to be selected by experience or estimated based on sufficient offline studies. For more adaptive model reduction, this paper presents a model reduction approach using the participation factors of state variables on the system modes, compared with the nonlinear model reduction method in [9]. In this new approach, the state variables (and so the associated generators) which are the most involved in a specific mode are reduced by a hybrid nonlinear model reduction method and the remaining generators will be linearized. Some advantages of using participation factors in the model reduction approach are: 1) it is dimensionless and so its application can be easily generalized. 2) it is independent of the particular choice of initial condition and also the characteristic of a disturbance that occurred in the system. 3) the error between the reduced-order model and the original detailed model is substantially lower compared to the use of fully linearized model reduction and rotor angle deviation-based model reduction approaches. Hence, the main contribution of this paper is to select participation factors associated with the system state variables as criteria to decide which generator models have more participation in dominant modes and their models should be considered nonlinear. The proposed approach is implemented on a real power system grid and the simulation results show the validity and benefit of the proposed approach.

In the rest of this paper, Section II briefly describes modal analysis in power system grid which includes the calculation of eigenvalues, eigenvectors and participation factors. In section

This work was supported in part by the ERC Program of the NSF and DOE under grant EEC-1041877 and in part by the NSF grant ECCS-1553863.

III, the model reduction process is explained. In Section IV, the simulation results of the proposed approach on a real power system grid are presented and compared with fully linearized model reduction and an adaptive model reduction using rotor angle criteria. Finally, the conclusion is explained in Section V.

## II. MODAL ANALYSIS

### A. State space representation

A power system dynamic model can be presented as a set of non-linear and first-order differential equations as the following form:

$$\dot{x} = f(x, u) \quad (1)$$

$$y = g(x, u) \quad (2)$$

$$x = \begin{bmatrix} x_1 \\ x_2 \\ \vdots \\ x_i \end{bmatrix}, \quad u = \begin{bmatrix} u_1 \\ u_2 \\ \vdots \\ u_l \end{bmatrix} \quad (3)$$

where vector $x$ is the state vector. $x_i$ and $u_l$ are states variables and inputs, respectively. The vectors $u$ and $y$ show the input and output vector of the system, respectively. $g$ is included nonlinear functions which relate state variables and inputs to the system outputs.

In the system equilibrium point, all the state variables are constant and their derivatives should be zero. Therefore, the following equation should be fulfilled.

$$f(x_0) = 0 \quad (4)$$

where $x_0$ is the state vector at the equilibrium point.

If the system deviates only by a very small amount from its equilibrium point, the nonlinear function can be linearized using the first terms of Taylor's series approximation. Hence, the linear formulation of (1) and (2) are:

$$\Delta \dot{x} = A\Delta x + B\Delta u \quad (5)$$

$$\Delta y = C\Delta x + D\Delta u \quad (6)$$

where matrices $A$, $B$, $C$ and $D$ are partial derivatives of functions $f$ and $g$ with respect to state vector $x$ and input vector $u$.

The state-space representation of the system in the frequency domain can be obtained as follows.

$$s\Delta x(s) - \Delta x(0) = \mathbf{A}\Delta x(s) + \mathbf{B}\Delta u(s) \quad (7)$$

$$\Delta y(s) = \mathbf{C}\Delta x(s) + \mathbf{D}\Delta u(s) \quad (8)$$

### B. Eigenvalues and Eigenvectors

In the state-space representation, matrix $\mathbf{A}$ is unique to the system for a given equilibrium point in (7), (8), but matrixes $\mathbf{B}$, $\mathbf{C}$, and $\mathbf{D}$ are dependent on both the equilibrium point and the choice of system inputs and outputs. The state matrix $\mathbf{A}$, and more specifically the eigenvalues of A, define the system around a selected equilibrium point. These eigenvalues satisfy the equation:

$$\det(\mathbf{A} - \lambda_i \mathbf{I}) = 0 \quad (9)$$

where $\lambda_i$ are eigenvalues of the matrix system $\mathbf{A}$.

The vectors $\mathbf{\Phi}_i$ and $\mathbf{\Psi}_i$ which satisfy the following equations are called right eigenvector and left eigenvector, respectively [12].

$$\mathbf{A}\mathbf{\Phi}_i = \lambda_i \mathbf{\Phi}_i \quad i = 1, 2, \cdots, n \quad (10)$$

$$\mathbf{\Psi}_i \mathbf{A} = \lambda_i \mathbf{\Psi}_i \quad i = 1, 2, \cdots, n \quad (11)$$

### C. Participation factor

Participation factors are used to figure out the relative interaction between state variables and modes. They examine the observability of a mode in a state variable as well as the state variable's contribution to the mode.

As (12) shows, the elements of both the right and left eigenvectors are used to calculate participation factors. The participation factors are dimensionless values and they are unaffected by the units used to measure state variables.

$$p_{ki} = \mathbf{\Phi}_{ki}\mathbf{\Psi}_{ik} \quad (12)$$

where $p_{ki}$ shows the relative contribution of the kth state variable in the $i$-th mode and in the reverse direction [12].

## III. MODEL REDUCTION PROCESS

### A. Partitioned power system network

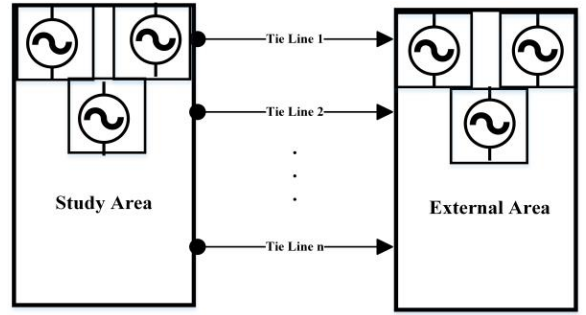

Figure 1. Partitioned power system model

As Fig.1 shows, the model reduction of a power system partitions the system model into two subsystems, including study area and external area. These two subsystems are connected through several tie lines. Each of these tie lines is considered as a constant voltage source generator and it injects electrical power from study area to external area and vice versa. As a result, the magnitude and phase angle of the output voltage at tie-line buses in each area are used as inputs to the adjacent subsystem.

In this paper, each generator and associated controllers are represented by totally nine first-order differential equations, including a detailed two-axis generator model, the first-order governor model, the non-reheat steam turbine model and the IEEE type-1 exciter.

When applying model reduction to the external area, it is required to specify the system's states and inputs. In this paper, the number of state variables and the inputs to the external area is defined as:

$$N_{state} = 9 N_{gen} \quad (13)$$

$$N_{in} = 2 N_{tie} \quad (14)$$

where $N_{state}$ is the number of state variables and $N_{in}$ is the number of inputs. $N_{gen}$ and $N_{tie}$ represent the number of generating units and the number of tie lines, respectively.

Nine differential equations of generating units can be formulated as nonlinear functions in (1) and (2). In this representation, the state vector $x$ and input vector $u$ are defined as:

$$x = (\delta \ P_m \ P_{gv} \ V_R \ R_f \ E_{fd} \ E_d^{'} \ E_q^{'} \ \omega)^T \quad (15)$$

$$u = (\theta \ V) \quad (16)$$

where $\delta$ and $\omega$ denote the rotor angle in rad and the speed of generators in rad/s, respectively. $P_m$ is the mechanical power, $P_{gv}$ is the governor output power, $V_R$ shows voltage regulator input and $R_f$ is rate feedback. $E_{fd}$, $E_q^{'}$, $E_d^{'}$ are field voltage, internal voltages on the q-axis and the d-axis, respectively. $\theta$ and $V$ are the voltage angle and voltage magnitude at boundary buses.

### B. Proposed Participation Factor-Based Adaptive Model Reduction

Given that the linearly reduced model performs satisfactorily under small disturbances, and typically a system is under small or no disturbance, it is reasonable to switch the type of model reduction for external area to keep a balance between the accuracy and speed of simulation [9]. The adaptive model reduction proposed in this paper improves the speed of simulation without loss of accuracy. In this approach, when the system is in the on-fault condition, the detailed original system is simulated; otherwise, the external area will be simplified. This process is similar to the adaptive nonlinear model reduction approach previously done in [9], but instead of rotor angle criteria, the modal analysis is done to decide which nonlinear functions in the external area could be reduced. To do this, eigenvalues, eigenvectors and participation factors of each state variable on the dominant system modes will be evaluated using equations (9) through (12). The level of interaction of each generator on the system modes is evaluated and they are arranged according to their absolute values. In this approach, one or several dominant modes energized by the fault will be identified. If the participation factor of the $k$th state variable into the $i$-th dominant mode, $p_{ki}$ in the external area reaches a predefined threshold $p_{max}$, the associated generator model is reduced by a hybrid nonlinear approach, while the generator models with small participation factors are fully linearized. This approach aims to maintain a compromise between accuracy and simulation speed following each disturbance. The algorithm of this method is shown in Fig.2. In this study, the balanced truncation method is used to obtain the reduced-order model of the external area.

The reduced model for generators with high participation factors using hybrid nonlinear model reduction approach can be expressed as follows:

$$\dot{\tilde{x}} = \mathbf{T} \begin{pmatrix} \tilde{f}(\tilde{\mathbf{T}}\tilde{x}, u) \\ \tilde{\mathbf{A}} \Delta \tilde{x} + \tilde{\mathbf{B}} \Delta u + \hat{x}_0 \end{pmatrix} \quad (17)$$

$$y = \tilde{\mathbf{T}} \tilde{x} \quad (18)$$

where $\tilde{f}$ includes nonlinear functions for generators with large participation factors in the external area and $\tilde{\mathbf{A}} \Delta \tilde{x} + \tilde{\mathbf{B}} \Delta u + \hat{x}_0$ is the representation of linearized functions for these generators. $\tilde{\mathbf{T}}$ is the inverse of transformation matrix $\mathbf{T}$. In this representation, $\hat{x}_0$ is the vector of the initial condition, $\tilde{\mathbf{A}}$ and $\tilde{\mathbf{B}}$ are defined as follows [9].

$$\tilde{\mathbf{A}} = \tilde{\mathbf{P}} \mathbf{A} \tilde{\mathbf{T}} \quad (19)$$

$$\tilde{\mathbf{B}} = \tilde{\mathbf{P}} \mathbf{B} \quad (20)$$

where $\tilde{\mathbf{P}}$ is reduced identity matrix.

The resulting model obtained by (17)-(20) has reduced nonlinearities, but it is not linear. In this approach, the nonlinear generator functions with large contribution to the study area are left nonlinear, while the nonlinear generator functions with small contribution are linearized. Therefore, the generators in the external area with a large contribution to the dynamics between the external area and study area are nonlinearly reduced.

### IV. SIMULATION RESULTS

The Northeastern Power Coordinating Council (NPCC) system with 48-machine, 140-bus system model [1] is studied as a test case. Fig. 3. shows the decomposition of the NPCC system into external and study areas. The study area is selected to be included 9 generators in the ISO-NE region with 81 state variables. The external area has 39 generators and 351 state variables.

The study area is preserved with nonlinear, fully detailed generation unit models, while the external area is respectively reduced using fully linearized, rotor angle deviation-based model reduction and participation factor-based adaptive methods for the purpose of comparison.

The fully linearized model and rotor angle deviation-based nonlinearly reduced models are generated from the previously described process in [9]. In the fully linearized model reduction approach, all generating units of the external area will be linearized as (5) and (6). The reduced-order models obtained by these three methods are compared with the original full-order model.

The simulations are performed in MATLAB R2021a on a computer with the Intel(R) Core (TM) i7-4790 CPU@ 3.60 GHz, 16.0 GB processor.

In this study, the total simulation time period is 16 seconds, the simulation time step is considered 0.01 seconds and a three-phase short circuit fault with duration of 0.39 seconds is applied to bus 3 of the NPCC system as shown in Fig. 3.

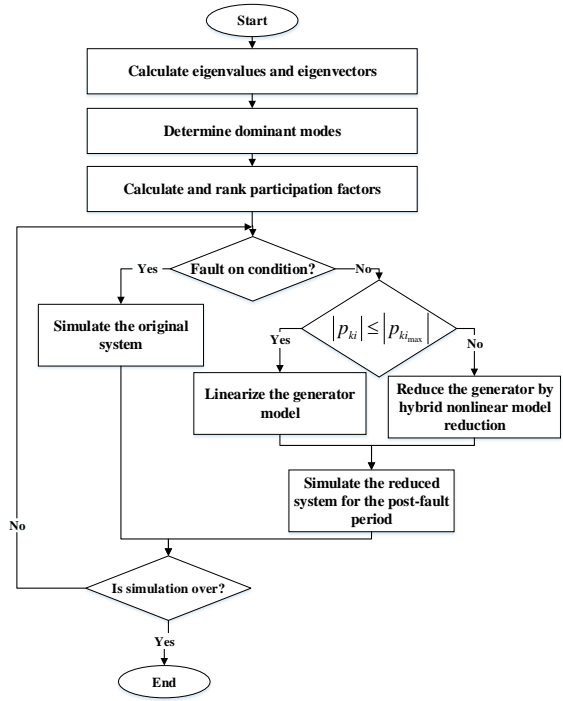

Figure 2. Proposed adaptive model reduction method on the external area

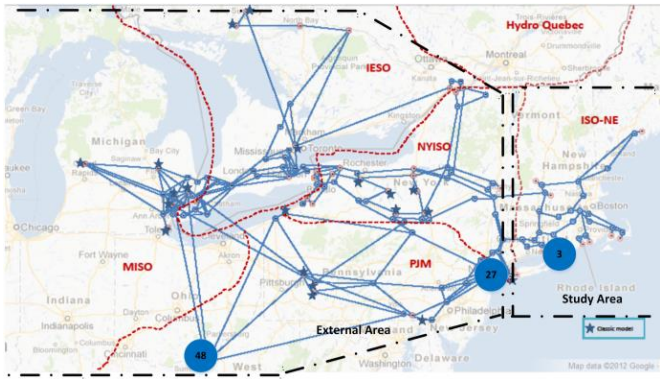

Figure 3. Partitioned NPCC system

The external area and its modal characteristics are important as far as they affect the power system analysis of the study area. Poorly damped oscillatory inter-area modes with low frequency (<1.0 HZ) oscillation and large magnitude are potentially dangerous and may have a detrimental impact on the dynamics of the study area. As a result, in this work, these modes are referred to as dominant modes.

The response of linearized power system model to zero input can be written as a linear combination of modes as follows [12]:

$$\Delta x(t) = \sum_{i=1}^{n} \mathbf{\Phi}_{ki} \mathbf{\Psi}_{ik} \Delta x(0)\, e^{\lambda_i t} \quad (21)$$

where $n$ denotes the number of modes and $\Delta x(0)$ is the initial condition. In this representation, the magnitude of mode $i$ can be calculated using (22). After calculation, two least damped modes with the largest magnitude are selected as the dominant modes in this paper.

$$Z_i = \mathbf{\Phi}_{ki} \mathbf{\Psi}_{ik} \Delta x(0) \quad (22)$$

The frequencies and damping ratios of these two modes are included in Table. I.

TABLE I. DOMINATED MODES'S FREQUENCY AND DAMPING RATIO

| Dominated Mode | Type of mode | Frequency (Hz) | Damping Ratio |
|---|---|---|---|
| Mode 1 | Inter-area mode | 0.6305 | 0.0804 |
| Mode 2 | Inter-area mode | 0.3874 | 0.1039 |

Based on experience, this paper considers a threshold of 0.5 to decide the effect of participation factor of each state variable on the dominant modes. Linearization is applied to nonlinear functions corresponding to generators with participation factors smaller than 0.5. This equates to the linearization of 37 generator models in the external area. Table II shows that generators 27 and 48 in the external area have the largest participation factors in two dominant modes and so their models should be nonlinearly reduced by (17) -(20). These two generators are marked in Fig.3.

TABLE II. PARTICIPATION FACTOR OF CRUCIAL GENERATORS

| Dominated Mode | PF of Crucial Generators | |
|---|---|---|
| | *Generator 27* | *Generator 48* |
| Mode 1 | 0.9978 | 0.0009 |
| Mode 2 | 0.5003 | 0.9996 |

Fig.4 demonstrates the rotor angle variations at bus 3 for generator 23 in the original full-order model, fully linearized model, rotor angle deviation-based adaptive reduced model and participation factor-based reduced model. It should be noticed that generator 23 is chosen to compare these three approaches because of its largest rotor angle deviations compared to the other generators.

From Fig. 4, the rotor angle mismatch error between the reduced-order model obtained by the fully linearized approach and original full-order model is relatively large, while the participation factor-based method is capable of closely following the rotor angle of the original full-order model.

In order to numerically compare the model reduction approaches, the root of mean squared errors for the state variables of generator 23 are determined using:

$$\varepsilon_i = \sqrt{\frac{\sum_{j=1}^{n}(x_{ij} - \hat{x}_{ij})^2}{n}} \quad (22)$$

where $x_{ij}$ is the *i*-th state variable's value at time step *j* and $\hat{x}_{ij}$ denotes the value of *i*-th state variable for the reduced-order model at time step *j*. *n* shows the required number of time steps.

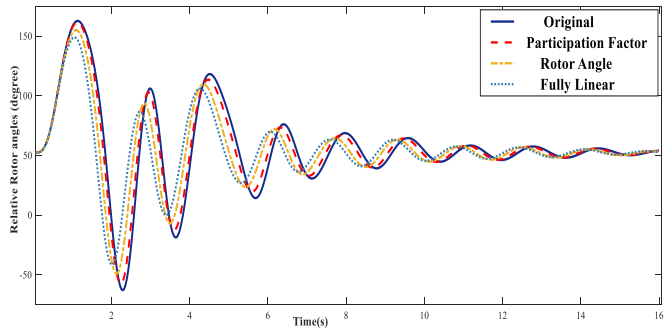

Figure 4. Rotor angle at bus 3 for generator 23

TABLE III. ROOT MEAN SQUARED ERROR OF STATE VARIABLES FOR GENERATOR 23

| States | Approach | | |
|---|---|---|---|
| | *Fully linear* | *Rotor-angle* | *Participation factor* |
| $\delta$, degrees | $2.59\times10^1$ | $17.13\times10^0$ | $5.77\times10^0$ |
| $P_m$, p.u. | $1.70\times10^{-3}$ | $1.70\times10^{-3}$ | $7.00\times10^{-4}$ |
| $P_{gv}$, p.u. | $1.98\times10^{-2}$ | $1.30\times10^{-2}$ | $4.50\times10^{-3}$ |
| $V_R$, p.u. | $1.71\times10^{-1}$ | $1.14\times10^{-1}$ | $4.02\times10^{-2}$ |
| $R_f$, p.u. | $1.34\times10^{-2}$ | $8.40\times10^{-3}$ | $3.10\times10^{-3}$ |
| $E_{fd}$, p.u. | $1.01\times10^{-1}$ | $6.50\times10^{-2}$ | $2.34\times10^{-2}$ |
| $E_d'$, p.u. | $7.09\times10^{-2}$ | $4.64\times10^{-2}$ | $1.61\times10^{-2}$ |
| $E_q'$, p.u. | $1.13\times10^{-2}$ | $7.20\times10^{-3}$ | $2.60\times10^{-3}$ |
| $\omega$, p.u. | $4.20\times10^{-3}$ | $2.80\times10^{-3}$ | $9.00\times10^{-4}$ |

Table III shows the calculated error $\varepsilon_i$ for each state variable. This result shows that for simulation of this contingency, the amount of error in the proposed participation factor-based method is significantly reduced compared to the fully linearized model. It also represents higher accuracy compared to the rotor angle deviation-based model reduction approach. Table IV compares the methodologies in terms of simulation time. The fully linearized approach is the fastest model with the lowest accuracy and the original full-order model is the most accurate model with the lowest speed. The purpose of adaptive model reduction is to maintain a compromise between accuracy and speed. The participation factor-based and rotor angle deviation-based model reduction approaches are in the same level of computation speed; however, the participation factor-based method suggests a higher degree of accuracy.

As already mentioned in [9], the threshold employed for the rotor angle deviation-based model reduction approach is system-specific and when the technique is implemented on a different system, the per-unit values for all criteria should be recalculated. Comparatively, the proposed model reduction approach utilizes participation factors, which are dimensionless values, so the approach can be more easily generalized.

TABLE IV. SIMULATION TIME FOR DIFFERENT APPROACHES

| System | Simulation time (Sec) |
|---|---|
| Original full-order | 1.092 |
| Linearly reduced | 0.44 |
| Partitioned rotor angle deviation-based | 0.52 |
| Partitioned participation factor-based | 0.56 |

## V. CONCLUSION

The proposed participation factor-based adaptive model reduction approach speeds up the power system transient stability simulation while retaining a reasonably high degree of accuracy. The method has been implemented on the partitioned NPCC 140-bus system and the results have been compared with a fully linearized model reduction approach and with a hybrid nonlinear model reduction approach based on rotor angle deviations. The results obtained from the reduced-order participation factor-based model are in close agreement with the original fully detailed model and the error in this method has substantially been decreased compared to the two other approaches.